\documentclass {tran-l}
 \newtheorem{thm}{Theorem}[section]
 \newtheorem{cor}[thm]{Corollary}
 \newtheorem{lem}[thm]{Lemma}
 \newtheorem{prop}[thm]{Proposition}
 \theoremstyle{definition}
 \newtheorem{defn}[thm]{Definition}
 \theoremstyle{remark}

\numberwithin{equation}{section}
\begin{document}
\title[Converting between Stirling, Tanh, Lah Numbers-Cumulants]
{Inversions relating Stirling, Tanh, Lah Numbers and an
application to Mathematical Statistics}

\author{Giacomo Della Riccia} 

\address{University of Udine, Department of Mathematics and Computer Science, Via delle Scienze 206. 33100-Udine,
Italy.} \email{dlrca@uniud.it}



\begin{abstract}
Inversion formulas have been found, converting between
\emph{Stirling}, \emph{tanh} and \emph{Lah} numbers. \emph{Tanh}
and \emph{Lah} polynomials, analogous to the Stirling
polynomials, have been defined and their basic properties
established. New identities for Stirling and tangent numbers and
polynomials have been derived from the general inverse relations.
In the second part of the paper, it has been shown that if
shifted-gamma probability densities and negative binomial
distributions are \emph{matched} by equating their first three
semi-invariants (cumulants), then the cumulants of the two
distributions are related by a pair of reciprocal linear
combinations equivalent to the inversion formulas established in
the first part.
\end{abstract}
\maketitle
\section{Introduction}
The usual form of an \emph{inversion formula} is
\begin{equation*}
g(n)=\sum_ia(n,i)f(i)\leftrightarrow f(n)=\sum_iA(n,i)g(i),
\end{equation*}
where $\left\{a(n,m), A(n,m)\right\}$ is a pair of \emph{inverse
numbers} and $\left\{f(n)\right\}$, $\left\{g(n)\right\}$ are
numerical sequences. An extensive study of inverse relations and
related topics can be found in \cite{Riorb}. In this paper we will
consider inversion formulas of a more general form:
\begin{equation*}
g(n,m)=\sum_ia(n,i)f(i,m)\leftrightarrow
f(n,m)=\sum_iA(n,i)g(i,m),
\end{equation*}
with \emph{double-sequences} $\left\{f(n,m)\right\}$,
$\left\{g(n,m)\right\}$ and corresponding \emph{number inverses}
$\left\{F(n,m)\right\}$, $\left\{G(n,m)\right\}$; in other words,
we are interested in inverse relations converting between
\emph{arrays}. Of course, from the general formula, one can get
inverse relations for \emph{sequences} by fixing one of the two
indices. As we shall see, the case $m=1$ is of particular
interest because, in general, it is associated with identities
involving important classical numbers. We used this approach to
find relations connecting Stirling $\left\{{n \brack m}, {n \brace
m}\right\}$, $\left\{\mbox{arctan}\ t(n,m), \mbox{tangent}\
T(n,m)\right\}$ and Lah $L(n,m)$ numbers. Then, in analogy to the
Stirling polynomials $\sigma_k(x)$ case, we introduced
\emph{tangent} $\delta_k(x)$ and Lah $\lambda_k(x)$ polynomials
and derived connecting relations. Most of the identities we
obtained seem to be original. Finally, we discussed the above
results in the light of a problem in Statistical Mathematics
dealing with semi-invariants (cumulants) of \emph{shifted-gamma}
probability densities $g(\vartheta;\ a,b,c)=\Gamma(\vartheta+c;\
a,b)$ and \emph{negative binomial distributions} $nb(\varpi;
r,\lambda)$. We showed that if $g(\vartheta;\ a,b,c)$ and
$nb(\varpi; r,\lambda)$ are \emph{matched} by equating their
first three cumulants, then the cumulants $\gamma(n)$ and
$\eta(n)$ of the two distributions are related by reciprocal
linear combinations equivalent to the \emph{array} inversion
formulas established previously.
\section{Relations between Stirling, tanh and Lah numbers}
We will use Stirling ${n\brack m}$, ${n\brace m}$ \cite{Kn},
arctan $t(n,m)$ and tangent $T(n,m)$ (\cite{Comt}, see pp.
258-260), and Lah $L(n,m)$ (\cite{Riora}, see pp. 43-44) numbers,
but with scale factors and appropriate notations:
\begin{eqnarray*}
&&\theta(n,m)=(-1)^{(n-m)/2}t(n,m),\quad
\Theta(n,m)=(-1)^{(n-m)/2}T(n,m);\\
&& \overline{{n\brack m}}=(-2)^{n-m}{n\brack m},\quad
\overline{{n\brace m}}=2^{n-m}{n\brace m};\\
&&\overline{l}(n,m)=(-1)^{n-m}\overline{L}(n,m)=(-1)^nL(n,m)=\frac{n!}{m!}{{n-1}\choose{m-1}}.
\end{eqnarray*}

The orthogonal relations satisfied by $\left\{\theta(n,m),
\Theta(n,m)\right\}$ and $\left\{{\overline{n\brack m}},
\overline{{n\brace m}}\right\}$ are of the form
$\sum_{i}a(n,i)A(i,m)=\sum_{i}A(n,i)a(i,m)=[m=n]$ and it is
easily verified by direct calculation that these relations are
also valid for $\left\{\overline{l}(n,m),
\overline{L}(n,m)\right\}$.

Recursions and egfs are simply obtained by introducing scales in
the ordinary relations. Since the egfs of $\theta(n,m)$,
$\Theta(n,m)$ egfs are of particular interest, we give the
following explicit derivation. The known egfs \cite{Comt} are
\begin{equation*}
t_m(v)=\sum_{n}t(n,m)\frac{v^n}{n!} =\frac{1}{m!}\arctan^m v,\quad
T_m(v)=\sum_{n}T(n,m)\frac{v^n}{n!}=\frac{1}{m!}\tan^m v.
\end{equation*}
Since $t_m(v)$, $T_m(v)$, as functions of $v$, and $m$ have
opposite parity, $t(n,m)=T(n,m)=0$ when $n-m$ is \emph{odd}. Thus,
with $n-m$ is even, putting $v=\imath u$, $\imath^2=-1$,
$\arctan(\imath u)=\imath\arg\tanh u$ and $\tan(\imath
u)=\imath\tanh u$, we get
\begin{eqnarray}\label{egf}
&\theta_m(u)=\sum_{n}\imath^{n-m}t(n,m)\frac{u^n}{n!}
=\sum_{n}(-1)^{\frac{n-m}{2}}t(n,m)\frac{u^n}{n!}=\frac{1}{m!}\arg\tanh^m u\\
&=\frac{1}{m!}\left(\frac{1}{2}\ln \frac{1+u}{1-u}\right)^m
=\frac{1}{m!}\left(\sum_{j}\theta_{2j+1}\frac{u^{2j+1}}{(2j+1)!}\right)^m;\ \theta_{2j+1}
=(-1)^{j}t_{2j+1}=(2j)!;\nonumber\\
&\Theta_m(u)_m(u)=\sum_{n}\imath^{n-m}T(n,m)\frac{u^n}{n!}
=\sum_{n}(-1)^{\frac{n-m}{2}}T(n,m)\frac{u^n}{n!}=\frac{1}{m!}\tanh^m
u\nonumber\\
&=\frac{1}{m!}\left(\sum_{j}\Theta_{2j+1}\frac{u^{2j+1}}{(2j+1)!}\right)^m;\
\Theta_{2j+1}=(-1)^{j}T_{2j+1}
=\frac{4^{j+1}(4^{j+1}-1)}{2j+2}B_{2j+2}.\nonumber
\end{eqnarray}
Thus, scale factors $(-1)^{(n-m)/2}$ change \emph{arctan}
$t(n,m)$, \emph{tangent} $T(n,m)$ numbers in \emph{arctanh}
$\theta(n,m)$, \emph{tanh} $\Theta(n,m)$ numbers. For brevity,
the pairs $\{\overline{{n\brack m}}, \overline{{n\brace m}}\}$,
$\{\theta(n,m), \Theta(n,m)\}$ and $\{\overline{l}(n,m),
\overline{L}(n,m)\}$ will be called \emph{Stirling}, \emph{tanh}
and \emph{Lah} numbers, respectively. Scaled numbers basic
properties are listed in Table 1.

\begin{table}[ht]
\begin{center}
\begin{tabular}{c}
\hline
$\mbox{Recurrence relations}$\\
$\overline{{{n+1}\brack m}}=\overline{{n\brack{
m-1}}}-2n\overline{{n\brack m}},\quad \overline{{{n+1}\brace
m}}=\overline{{n\brace{ m-1}}}+2m\overline{{n\brace m}};$\\
$\overline{{0\brack m}}=\overline{{0\brace m}}=[m=0],\quad
\overline{{n\brack 0}}=\overline{{n\brace 0}}=[n=0].$\\
$\theta(n+1,m)=\theta(n,m-1)+n(n-1)\theta(n-1,m),$\\
$\Theta(n+1,m)=\Theta(n,m-1)-m(m+1)\Theta(n,m+1);$\\
$\theta(0,m)=\Theta(0,m)=[m=0];\ \theta(n,0)=\Theta(n,0)=[n=0];\
\theta(1,0)=\Theta(1,0)
=0.$\\
$\overline{l}(n+1,m)=(n+m)\overline{l}(n,m)+\overline{l}(n,m-1),$\\
$\overline{L}(n+1,m)=-(n+m)\overline{L}(n,m)+\overline{L}(n,m-1);$\\
$\overline{l}(n,0)=\overline{L}(n,0)\equiv0;\quad
(-1)^m\overline{l}(0,m)=\overline{L}(0,m)=-\frac{1}{m!}+[m=0].$\\
\hline
$\mbox{Duality laws and Orthogonal relations}$\\
$\overline{{n\brack m}}=(-1)^{n-m}\overline{{-m\brace {-n}}};\
\theta(n,m)=\Theta(-m,-n);\
\overline{l}(n,m) =(-1)^{n-m}\overline{L}(-m,-n).$\\
$\sum_{i}a(n,i)A(i,m)=\sum_{i}A(n,i)a(i,m)=[m=n].$\\
\hline
$\mbox{Exponential generating functions}$\\
$\sum_{n}\overline{{n\brack m}}\frac{u^n}{n!}
=\frac{1}{m!}\left(\frac{1}{2}\ln (1+2u)\right)^m,\quad
\sum_{n}\overline{{n\brace m}}\frac{u^n}{n!}
=\frac{1}{m!}\left(\frac{e^{2u}-1}{2}\right)^m$\\
$\theta_m(u)=\sum_{n}\theta(n,m)\frac{u^n}{n!}
=\frac{1}{m!}\arg\tanh^m u=\frac{1}{m!}\left(\frac{1}{2}\ln
\frac{1+u}{1-u}\right)^m$\\
$\Theta_m(u)=\sum_{n}\Theta(n,m)\frac{u^n}{n!}=\frac{1}{m!}\tanh^m
u$\\
$\sum_{n}\overline{l}(n,m)\frac{u^n}{n!}
=\frac{1}{m!}\left(\frac{u}{1-u}\right)^m,\quad
\sum_{n}\overline{L}(n,m)\frac{u^n}{n!}
=\frac{1}{m!}\left(\frac{u}{1+u}\right)^m.$\\
\hline
\end{tabular}
\end{center}
\caption{Basic properties of scaled Stirling, tanh and Lah
numbers}
\end{table}
\begin{thm}\label{conv}
Numbers in each pair, \emph{Stirling}, \emph{tanh} and \emph{Lah},
convert between numbers in the other two pairs.
\end{thm}
Let $\tanh u=\frac{v}{1+v}$, $v=\frac{e^{2u-1}}{2}$, then from the
egfs listed in Table 1 we get:
\begin{equation}
\begin{array}{l}\label{sequence}
\sum_{n}\Theta(n,m)\frac{u^n}{n!}=\frac{1}{m!}\tanh^m u
=\frac{1}{m!}\left(\frac{v}{1+v}\right)^m=\sum_i\overline{L}(i,m)\frac{v^i}{i!}\\
=\sum_i\overline{L}(i,m)\frac{1}{i!}\left(\frac{e^{2u}-1}{2}\right)^i=
\sum_i\overline{L}(i,m)\sum_n\overline{{n\brace
i}}\frac{u^n}{n!}\\
=\sum_n\left(\sum_i\overline{L}(i,m)\overline{{n\brace
i}}\right)\frac{u^n}{n!}.
\end{array}
\end{equation}
These equations imply
$\Theta(n,m)=\sum_{i=m}^n\overline{L}(i,m)\overline{{n\brace i}}$
and, by dualities and inversions (Table 1),
$\theta(n,m)=\sum_{i=m}^n\overline{l}(i,m)\overline{{n\brack
i}}$, $\overline{{n\brack m}}=\sum_i\overline{L}(n,i)\theta(i,m)$
and $\overline{{n\brace m}}=\sum_i\overline{l}(n,i)\Theta(i,m)$.
Hence, Lah numbers convert between Stirling and tanh numbers. In
Table 2 we listed the identities that are derived by use of
inversions and/or dualities given in Table 1. Since in the third
and fourth inversion formulas Stirling numbers convert between
Lah and tanh numbers, and in the fifth and sixth tanh numbers
convert between Lah and Stirling numbers, the proof is complete.

The basic structure connecting tanh and Stirling numbers is the
following.
\begin{cor}\label{coeff}
Tanh numbers are finite sums of multiples of Stirling numbers, and
inversely
\begin{equation*}
\begin{array}{l}
\theta(n,n-k)=\sum_{i=0}^ki!{{n}\choose i}{{n-1}\choose
i}\overline{{{n-i}\brack {n-k}}},\\
\overline{{n\brack {n-k}}}=\sum_{i=0}^k(-1)^ii!{{n}\choose
i}{{n-1}\choose i}\theta(n-i,n-k),\\
\Theta(n+k,n)=\sum_{i=0}^k(-1)^{k-i}i!{{n+i}\choose{i}}
{{n+i-1}\choose{i}}\overline{{{n+k}\brace{n+i}}},\\
\overline{{n+k\brace
n}}=\sum_{i=0}^k(-1)^{k-i}i!{{n+i}\choose{i}}{{n+i-1}\choose{i}}\Theta(n+k,n+i).
\end{array}
\end{equation*}
\end{cor}
These relations are obtained from the first two inverse pairs in
Table 2 with $k=n-m$, $i$ replaced by $n-i$ and the use of Lah
numbers explicit expressions.
\begin{table}[ht]
\begin{center}
\begin{tabular}{|c|}
\hline $\overline{{n\brack m}}=(-2)^{n-m}{n\brack m},\quad
\overline{{n\brace
m}}=2^{n-m}{n\brace m}$\\
$\theta(n,m)=(-1)^{(n-m)/2}t(n,m),
\quad \Theta(n,m)=(-1)^{(n-m)/2}T(n,m)$\\
$\overline{l}(n,m)=(-1)^nL(n,m),\quad \overline{L}(n,m)
=(-1)^mL(n,m)$\\
\hline $\Theta(n,m)=\sum_i\overline{L}(i,m)\overline{{n\brace
i}}\quad \leftrightarrow\quad \overline{{n\brace
m}}=\sum_i\overline{l}(i,m)\Theta(n,i)$\\$\theta(n,m)=\sum_i\overline{l}(n,i)\overline{{i\brack
m}}\quad \leftrightarrow\quad \overline{{n\brack
m}}=\sum_i\overline{L}(n,i)\theta(i,m)$\\
\hline $\overline{L}(n,m) =\sum_i\overline{{n\brack
i}}\Theta(i,m)\quad \leftrightarrow\quad
\Theta(n,m)=\sum_i\overline{{n\brace
i}}\overline{L}(i,m)$\\
$\overline{l}(n,m)=\sum_i\overline{{i\brace m}}\theta(n,i)\quad
\leftrightarrow\quad \theta(n,m)=\sum_i\overline{{i\brack
m}}\overline{l}(n,i)$\\
\hline $\overline{l}(n,m)=\sum_i\theta(n,i)\overline{{i\brace
m}}\quad \leftrightarrow\quad \overline{{n\brace m}}=\sum_i\Theta(n,i)\overline{l}(i,m)$\\
$\overline{L}(n,m) =\sum_i\Theta(i,m)\overline{{n\brack i}}\quad
\leftrightarrow\quad \overline{{n\brack m}}=\sum_i\theta(i,m)\overline{L}(n,i)$\\
\hline
\end{tabular}
\end{center}
\caption{Conversions between Stirling, tanh, and Lah numbers}
\end{table}
As we know (\cite{Knu}, see p. 418), Stirling numbers
$\overline{{x\brack {x-k}}}$, $\overline{{x+k\brace x}}$ can be
viewed as polynomial in $x$. Thus, Corollary (\ref{coeff}) implies
that $\theta(x,x-k)$, $\Theta(x+k,x)$ can also be treated as
polynomials; these have the following properties.
\begin{prop}\label{poly1}
If\ $k=2j\geq0$, then $\theta(x,x-k)$, $\Theta(x+k,x)$ are
polynomials in $x$ having degree $\frac{3k}{2}=3j$ and leading
coefficient $\frac{1}{3^j\times j!}$, $(-1)^j\frac{1}{3^j\times
j!}$, respectively.
\end{prop}
The proof is by induction on $k$ applied to tanh numbers
recurrence relations (Table 1) written with $k=n-m$ notations.
Details of the proof are omitted because they are the same as
those used by Gessel and Stanley (\cite{Gess}, see p. 25 ) in
their study on Stirling numbers structure.

The first few cases are the following.
\begin{eqnarray*}
&\theta(x,x-2)=-\Theta(x,x-2)=\frac{3!}{3\times1!}{x\choose 3}\\
&\theta(x,x-4)=\frac{6!}{3^2\times2!}{{x+1}\choose 6}-2^4{x\choose
5};\quad \Theta(x,x-4)=\frac{6!}{3^2\times2!}{{x+1}\choose
6}-2^3\times3{x\choose 5}.
\end{eqnarray*}
As pointed out in the Introduction, the general inverse relations
in Table 2 yield interesting results in the case of $m=1$,
essentially because
\begin{eqnarray*}
&\overline{{n\brack 1}}=(-2)^{n-m}(n-1)!;\quad \overline{{n\brace
1}}=2^{n-m}.\\
&\theta(2n+1,1)=(2n)!;\quad \Theta(2n+1,1)=(-1)^nT_{2n+1};\quad
\theta(2n,1)=\Theta(2n,1)\equiv0.
\end{eqnarray*}
The first general pair of inverse relations gives Stirling
numbers identities:
\begin{equation}\label{nids}
\begin{array}{c}
\sum_{i=1}^{n}(-1)^{i-1}i!\overline{{{n}\brace{i}}}=
\begin{cases}
0,\ \mbox{even}\ n\\
\Theta(n,1),\ \mbox{odd}\ n\\
\end{cases}
\end{array}
\leftrightarrow\quad 2^{n-1}=\sum_{i=1}^n i!\Theta(n,i).
\end{equation}
The identity when $n$ is even was given for the first time by
Lengyel (\cite {Len}, see p. 7), whereas the identity for $n$ odd
is new. The third inverse pair yields:
\begin{equation}\label{nidp}
2^{n-1}=\sum_{i=1}^n\Theta(n,i)i!\quad \leftrightarrow\quad
n!=\sum_{i=1}^{n}\theta(n,i)2^{i-1},
\end{equation}
showing that tanh numbers convert between factorials and powers of
2. From the fourth we get an inverse pair:
\begin{equation*}
\theta(n,1)
=\sum_{i=1}^n\overline{l}(n,i)(-2)^{i-1}(i-1)!\leftrightarrow
(-2)^{n-1}(n-1)! =\sum_{i=1}^n\overline{L}(n,i)\theta(i,1),
\end{equation*}
extending an inversion formula
$n!=\sum_{i=1}^n\overline{L}(n,i)2^{i-1}i!\leftrightarrow
2^{n-1}n!=\sum_{i=1}^n\overline{l}(n,i)i!$ given by Lah
(\cite{Lah}, see p. 207). Finally, the fifth and the sixth pairs
disclose original identities involving two out of the Stirling,
tanh and Lah numbers:
\begin{eqnarray*}
&\Theta(n,1)=\sum_{i=1}^n\overline{{{n}\brace{i}}}(-1)^{i-1}i!\quad
\leftrightarrow\quad (-1)^{n-1}n!=\sum_{i=1}^n
\overline{{{n}\brack{i}}}\Theta(i,1)\\
&(-1)^{n-1}n!=\sum_{i=1}^n\Theta(i,1)\overline{{{n}\brack{i}}}\
\leftrightarrow\ (-2)^{n-1}(n-1)!=
\sum_{i=1}^n\theta(i,1)\overline{L}(n,i).
\end{eqnarray*}
\section{Stirling, tanh and Lah polynomials}
Since $\theta(x,x-k)=\sum_{i=0}^ki!{{x}\choose i}{{x-1}\choose
i}\overline{{{x-i}\brack {x-k}}}$ and
$\overline{l}(x,x-k)=k!{{x}\choose{k}}{{x-1}\choose{k}}$ vanish
for $x=0,1,\ldots,k$ , \emph{tanh polynomials} $\delta_k(x)$ and
\emph{Lah polynomials} $\lambda_k(x)$ can be defined by rules
similar to $\overline{\sigma}_k(x)=\left.\overline{{x\brack
{x-k}}}\right/x^{\underline{k+1}}$ \cite{Kn} used for Stirling
polynomials. For clarity, we recall that
$x^{\underline{k+1}}=x(x-1)\ldots,(x-k)=(-1)^{k+1}(k-x)^{\underline{k+1}}$.
\begin{defn}\label{stp1}
\begin{eqnarray*}
&\delta_k(x)=\frac{\theta(x,x-k)}{x^{\underline{k+1}}}\
\sim\ \delta_k(k-x)=(-1)^{k+1}\frac{\Theta(x,x-k)}{x^{\underline{k+1}}};\\
&\lambda_k(x)=\frac{\overline{l}(x,x-k)}{x^{\underline{k+1}}}\
\sim\
\lambda_k(k-x)=-\frac{\overline{L}(x,x-k)}{x^{\underline{k+1}}}=(-1)^{k+1}\lambda_k(x).
\end{eqnarray*}
\end{defn}
In the case of integers $n,m$, the above definitions assume the
form
\begin{equation}\label{stp2}
\theta(n,m)=\frac{n!}{(m-1)!}\delta_{n-m}(n);\quad
\Theta(n,m)=-\frac{n!}{(m-1)!}\delta_{n-m}(-m).
\end{equation}
Corollary (\ref{coeff}) with $x$ instead of $n$ and a factor
$x^{\underline{k+1}}$ divided out, yields:
\begin{prop}\label{coeff3}
Tanh polynomials $\delta_k(x)$ are finite sums of multiples of
Stirling polynomials $\overline{\sigma}_k(x)$, and inversely
\begin{eqnarray*}
&\delta_k(x)=\sum_{i=0}^k{{x-1}\choose
{i}}\overline{\sigma}_{k-i}(x-i)\ \sim\
\delta_k(k-x)=\sum_{i=0}^k{{k-x-1}\choose
{k-i}}\overline{\sigma}_i(i-x);\\
&\overline{\sigma}_k(x)=\sum_{i=0}^k{{x-1}\choose
{i}}\delta_{k-i}(x-i)\ \sim\
\overline{\sigma}_k(k-x)=\sum_{i=0}^k{{k-x-1}\choose
{k-i}}\delta_i(i-x).
\end{eqnarray*}
\end{prop}
\begin{prop}
The generating functions of $x\delta_k(x)$ and $x\delta_{k}(k+x)$
are:
\begin{eqnarray*}
&\sum_{k}x\delta_k(x)u^{k}=\left(u\coth
u\right)^x=\left(\sum_{j}2^{2j}\frac{B_{2j}}{(2j)!}u^{2j}\right)^x,\
B_{2j}\ \emph{Bernoulli numbers};\\
&\sum_{k}x\delta_{k}(k+x)u^{k}=\left(\frac{1}{u}\arg\tanh
u\right)^x=\left(\frac{1}{2u}\ln
\frac{1+u}{1-u}\right)^{x}=\left(\sum_{j\geq0}\frac{1}{2j+1}u^{2j}\right)^{x}.
\end{eqnarray*}
\end{prop}
The egf of $\theta(n+k,n)$ (\ref{egf}) and definitions
(\ref{stp2}) yield
\begin{eqnarray*}
&&\sum_{k}\theta(n+k,n)\frac{u^{n+k}}{(n+k)!}
=\sum_{k}\frac{(n+k)!}{(n-1)!}\delta_k(n+k)\frac{u^{n+k}}{(n+k)!}=\frac{1}{n!}\arg\tanh^n
u,\\
&&\sum_{k}x\delta_{k}(k+x)u^{k}=\left(\frac{1}{u}\arg\tanh
u\right)^x
\end{eqnarray*}
where, on the basis of the "polynomial" argument, we replaced $n$
with $x$. Using the egf of $\Theta(n+k,n)$ and proceeding as
above, we obtain $\sum_{k}x\delta_k(x)u^{k}=\left(u\coth
u\right)^x$.
\begin{prop}\label{poly2}
Tanh polynomials $\delta_{k}(x)$ satisfy the recurrence relation
\begin{eqnarray*}
&(x+1)\delta_{k}(x+1)=(x-k)\delta_{k}(x)+(x-1)\delta_{k-2}(x-1);\quad
x\delta_0(x)\equiv1.\\
&\delta_{k}(x),\ k=2j>0,\ \emph{has degree}\ k/2-1\ \emph{and}\
\delta_{k}(x)\equiv0,\ k=2j+1.
\end{eqnarray*}
\end{prop}
The recurrence relation is obtained by dividing out a common
factor $x^{\underline{k}}$ in the recurrence relation of
$\theta(n,m)$ (Table 1), written with $n=x$ and $m=x+1-k$.
(Compare with
$(x+1)\sigma_k(x+1)=(x-k)\sigma_k(x)+x\sigma_{k-1}(x)$,
(\cite{Kn} Exercise 18, Chapter 6)). $\delta_k(x),\ k=2j>0 $ has
degree $k/2-1$ follows from Proposition (\ref{poly1}) and
$\delta_{k}(x)\equiv0,\ k=2j+1$, because
$\theta(x,x-k)=\Theta(x+k,x)\equiv0,\ k=2j+1$. (Recall that
$\overline{\sigma}_k(x),\ k>0$, has degree $k-1$,\ \cite{Kn}). The
first few cases are the following:
\begin{eqnarray*}
&\delta_{k}(1)=2^{k}\frac{B_{k}}{k!}+[k=1];\quad k\delta_k(0)
=2^k(2^k-2)\frac{B_k}{k!},\quad k>0;\\
&\delta_{k}(-1)=-\frac{\Theta_{k+1}}{(k+1)!};\quad
\delta_{2j}(2j+1)=\frac{1}{2j+1}.\\
&\delta_2(x)\equiv\frac{1}{3};\ \delta_4(x) =\frac{1!}{3^2\times
2!}{{x-1}\choose 1}-\frac{1}{3^2\times5};\
\delta_6(x)=\frac{2!}{3^3\times3!}{x\choose
2}-\frac{2^3}{3^4\times5}{{x-1}\choose
1}+\frac{2}{3^3\times7\times5}.
\end{eqnarray*}
A companion expression of (\ref{nids}), derived from Proposition
(\ref{coeff3}), is
\begin{equation*}
\begin{array}{c} \sum_{i=1}^{k}{{x-1}\choose
{i}}\overline{\sigma}_{k-i}(x-i)=
\begin{cases}
0,\quad \mbox{odd}\ k\\
\delta_k(x),\quad \mbox{even}\ k\\
\end{cases}
\end{array};\quad \overline{\sigma}_k(x)=\sum_{i=1}^{k}{{x-1}\choose
{i}}\delta_{k-i}(x-i).
\end{equation*}
$\sum_{i=1}^{2j+1}{{x-1}\choose
{i}}\overline{\sigma}_{2j+1-i}(x-i)=0$ is a new Stirling
polynomials identity.

The properties of $x\lambda_k(x)$ follow at ounce from
$x\lambda_k(x)={x\choose k}$.
\begin{prop}
Polynomials $x\lambda_k(x)={x\choose k}$ have degree $k$. Their
generating function and recurrence relation are
\begin{equation*}
\begin{array}{l}
\sum_{k\geq0}x\lambda_k(x)u^{k}=(1+u)^x,\quad
\sum_{k\geq0}x\lambda_{k}(k+x)u^{k}=\frac{1}{(1-u)^x},\\
(x+1)\lambda_k(x+1)={{x+1}\choose k}={{x}\choose
{k-1}}+{{x}\choose
{k}}=x\left[\lambda_{k-1}(x)+\lambda_k(x)\right],\quad
x\lambda_0(x)\equiv1.
\end{array}
\end{equation*}
\end{prop}
\section{An application to a problem in Mathematical Statistics}
We now show that the above results have an application in a
problem of Statistical Mathematics dealing with semi-invariants
(cumulants) of \emph{shifted-gamma} densities $g(\vartheta;\
a,b,c)=\Gamma(\vartheta+c;\ a,b)$ and \emph{negative binomial
distributions} $nb(\varpi; r,\lambda)$.

For the sake of completeness, we first recall some standard
definitions.
\begin{eqnarray*}
&\begin{array}{l} g(\vartheta;\ a,b,c)=\Gamma(\vartheta+c;\
a,b)=\begin{cases} \frac{1}{b^a\Gamma(a)}(\vartheta+c)^{a-1}
\exp[-\frac{(\vartheta+c)}{b}],&\text{$\vartheta>-c$},\\
    0&\text{otherwise}
  \end{cases}\\
  \end{array}\\
&nb(\varpi;\
r,\lambda)=\frac{\Gamma(r+\varpi)}{\Gamma(r)\Gamma(\varpi+1)}
\left(\frac{1}{1+\lambda}\right)^r\left(\frac{\lambda}{1+\lambda}\right)^{\varpi},\quad
\lambda\geq0
\end{eqnarray*}
From the \emph{moment} egfs
$M_{sg}(t)=e^{-ct}\left/(1-bt)^a\right.$,
$M_{nb}(t)=1\left/[1-\lambda(e^t-1)]^r\right.$ and the
\emph{cumulant} egfs $\ln M_{sg}(t)$, $\ln M_{nb}(t)$:
\begin{equation*}
\begin{array}{l}
\ln
M_{sg}(t)=\sum_{n>0}\gamma(n)\frac{t^n}{n!}=(-c+ab)t+\sum_{n>1}a\frac{b^n}{n}t^n,\\
\ln M_{nb}(t)=\sum_{n>0}\eta(n)\frac{t^n}{n!}
=r\sum_{m>0}\frac{\lambda^m}{m}(e^t-1)^m\\
=\sum_{n>0}\left[\sum_{m>0}r(m-1)!{n\brace
m}\lambda^m\right]\frac{t^n}{n!}\ \left[\mbox{with}\
(e^t-1)^m=m!\sum_{n\geq0}{n\brace m}\frac{t^n}{n!}\right],
\end{array}
\end{equation*}
we get the $n$th cumulants $\gamma(n)$, $\eta(n)$ of
$g(\vartheta;\ a,b,c)$, $nb(\varpi;\ r,\lambda)$:
\begin{eqnarray}
&\gamma(n)=
\begin{cases}\label{c1}
    -c+ab, & \text{$n=1$},\\
    (n-1)!\ ab^n, & \text{$n>1$};
  \end{cases}\\
&\eta(n)=\sum_{m=1}^nr(m-1)!2^{m-n}\overline{{n\brace
m}}\lambda^{m},\quad n>0.\label{c2}
\end{eqnarray}
\begin{defn}Two distributions $g(\vartheta;\ a,b,c)$ and $nb(\varpi;\ r,\lambda)$
are said to be \emph{matched} if their first three cumulants are
equal.
\end{defn}
By equating the first three cumulants
\begin{equation*}
\begin{array}{l}
1st\ cumulant\ (\mbox{mean})=\ \mu=-c+ab=r\lambda\\
2nd\ cumulant\ (\mbox{variance})=\sigma^2=ab^2=r\lambda(1+\lambda)\\
3rd\ cumulant=\gamma(3)=2!\ ab^3=\eta(3)=2!\
r\lambda(1+\lambda)(1/2+\lambda)
\end{array}
\end{equation*}
we get matching conditions
\begin{equation}\label{match}
\begin{array}{lll}
a=r\frac{\lambda(1+\lambda)}{(1/2+\lambda)^2},
& b=1/2+\lambda, & c=\frac{ab}{1+2b}=\frac{r\lambda}{1+2\lambda};\\
r=\frac{ab^2}{b^2-1/4}, & \lambda=b-1/2.
\end{array}
\end{equation}
For convenience, we will use scaled cumulants
$\overline{\eta(n)}=\frac{2^n}{r}\eta(n)$ and
$\overline{\gamma(n)}=\frac{2^n}{r}\gamma(n)$.
\begin{lem}\label{le2}
In a matched pair $\left\{g(\vartheta;\ a,b,c), nb(\varpi;\
r,\lambda)\right\}$, cumulants $\overline{\gamma}(n)$ and
$\overline{\eta}(n)$ are polynomials in $\lambda$ having degree
$n$
\begin{eqnarray}
&\overline{\gamma}(n)=\sum_{m=1}^n\gamma(n,m)\lambda^{m},\quad
n>0,\quad \overline{\eta}(n) =\sum_{m=1}^n\eta(n,m)\lambda^{m};\nonumber\\
&\gamma(n,m)=2^m(m-1)!\left[\overline{l}(n-1,m-1)
+2m\overline{l}(n-1,m)\right],\label{q}\\
&\eta(n,m)=2^{m}(m-1)!\overline{{n\brace m}}.\label{c3}
\end{eqnarray}
\end{lem}
The lemma holds for $\overline{\eta}(n)$ as we already know
(\ref{c2}). It is also true for $\overline{\gamma}(n)$ because if
$\left\{a,b,c\right\}$ are replaced by the matching values
(\ref{match}), then from (\ref{c1}):
\begin{eqnarray*}
&\begin{array}{l}
\overline{\gamma}(n)=\begin{cases} 2\lambda,\quad n=1 & \\
2^n(n-1)!\lambda(1+\lambda)(1/2+\lambda)^{n-2}=\sum_{m=1}^n\gamma(n,m)\lambda^m
, & \text{$n>1$},
\end{cases}\\
\end{array}\\
&\lambda(1+\lambda)(1/2+\lambda)^{n-2}=\sum_{m=1}^n\left[{{n-2}\choose
{n-m}}\frac{1}{2^{n-m}} +{{n-2}\choose
{n-m-1}}\frac{1}{2^{n-m-1}}\right]\lambda^n\\
&\mbox{and}\ \gamma(n,m)=2^m(m-1)!\left[\overline{l}(n-1,m-1)
+2m\overline{l}(n-1,m)\right],\quad n\geq1.
\end{eqnarray*}
\begin{thm}\label{th2}
In a matched pair $\left\{g(\vartheta;\ a,b,c), nb(\varpi;\
r,\lambda)\right\}$, cumulants $\overline{\gamma}(n)$ and
$\overline{\eta}(n)$ of the two distributions are related by
reciprocal linear combinations
\begin{eqnarray*}
&&\overline{\gamma}(n+1)=\sum_{i=0}^n\theta(n,i)\overline{\eta}(i+1)\quad
\leftrightarrow\quad \overline{\eta}(n+1)
=\sum_{i=0}^n\Theta(n,i)\overline{\gamma}(i+1),\quad n\geq0,\\
&&\emph{where}\ \theta(n,m)=(-1)^{(n-m)/2}t(n,m)\ \emph{and}\
\Theta(n,m)=(-1)^{(n-m)/2}T(n,m).
\end{eqnarray*}
\end{thm}
Since cumulants are polynomials in $\lambda$, a solution
$\theta(n,m)$ and $\Theta(n,m)$ exists if like powers of
$\lambda$ on both sides are equal:
\begin{eqnarray}
&&\eta(n+1,m+1)\equiv[\lambda^{m+1}]\sum_i\Theta(n,i)\overline{\gamma}(i+1)
=\sum_i\Theta(n,i)\gamma(i+1,m+1),\label{sys1}\\
&&\gamma(n+1,m+1)\equiv[\lambda^{m+1}]\sum_i\theta(n,i)\overline{\eta}(i+1)
=\sum_i\theta(n,i)\eta(i+1,m+1).\nonumber
\end{eqnarray}
From (\ref{c3}) and Stirling recurrence relation, and (\ref{q}),
we get
\begin{equation}\label{cum}
\begin{array}{l}
\eta(n+1,m+1)=2^{m+1}m!\overline{{{n+1}\brace
{m+1}}}=2^{m+1}m!\left[\overline{{n\brace{
m}}}+2(m+1)\overline{{n\brace
{m+1}}}\right],\\
\gamma(i+1,m+1)=2^{m+1}m!\left[\overline{l}(i,m)
+2(m+1)\overline{l}(i,m+1)\right].
\end{array}
\end{equation}
Substituting (\ref{cum}) into (\ref{sys1}), rearranging the terms
and dividing by $2^{m}m!$, we obtain a recursion in $m$ which
automatically terminates after $n-m$ recursive steps
\begin{eqnarray*}
&\overline{{n\brace m}}-\sum_{i=m}^{n}\Theta(n,i)
\overline{l}(i,m)=2(m+1)\left[\overline{{n\brace
{m+1}}}-\sum_{i=m+1}^{n}\Theta(n,i)\overline{l}(i,m+1)\right]\\
&=2^2(m+2)(m+1)\left[\overline{{n\brace
{m+2}}}-\sum_{i=m+1}^{n}\Theta(n,i)\overline{l}(i,m+2)\right]=,\ldots,=0,
\end{eqnarray*}
Thus $\overline{{n\brace
m}}-\sum_i\Theta(n,i)\overline{l}(i,m)=0$, and
$\Theta(n,m)=(-1)^{(n-m)/2}T(n,m)$ from Theorem (\ref{conv}).
Since $\gamma(n+1)=\sum_{i=0}^n\theta(n,i)\eta(i+1)$ is the
reciprocal of $\overline{\eta}(n+1)
=\sum_{i=0}^n\Theta(n,i)\overline{\gamma}(i+1)$, $\theta(n,m)$
and $\Theta(n,m)$ are necessarily \emph{number inverses},
therefore $\theta(n,m)=(-1)^{(n-m)/2}t(n,m)$. This completes the
proof. Conversely,
\begin{equation*}
\eta(n+1,m+1)
=\sum_i\Theta(n,i)\gamma(i+1,m+1)\leftrightarrow\gamma(n+1,m+1)=\sum_i\theta(n,i)\eta(i+1,m+1)
\end{equation*}
implies Theorem (\ref{conv}) and the inversion relations in Table
2, hence, \emph{the problem in Theorem \emph{(\ref{th2})} of
converting between cumulants is equivalent to the problem in
Theorem \emph{(\ref{conv})} of converting between Stirling, tanh
and Lah numbers}.

Finally, let us show that the inverse pair (\ref{nidp}):
$n!=\sum_{i=0}^n\theta(n,i)2^{i-1}\ \leftrightarrow\
2^{n-1}=\sum_{i=0}^n\Theta(n,i)i!,$ corresponds to an interesting
particular case of Theorem (\ref{th2}). In fact, let
$\lambda\rightarrow0$ while $r\lambda=\alpha$ remains constant,
then, as we know, the negative binomial distribution
$nb(\varpi;r,\lambda)$ tends towards the Poisson distribution
$p(\varpi; \alpha)$ with cumulants all equal to $\alpha$. The
matched shifted-gamma density tends towards
$g(\vartheta;4\alpha,\frac{1}{2},\alpha)$, since, according to
the matching conditions (\ref{match}), $a\rightarrow 4\alpha,\
b\rightarrow\frac{1}{2},\ c\rightarrow\alpha$, and its cumulants
are obtained from (\ref{c1}). Using the cumulant limiting values,
one verifies that the reciprocal relations in Theorem (\ref{th2})
reduce to (\ref{nidp}).
\section{Main results and Conclusion}
From general inverse relations converting between Stirling, tanh
and Lah numbers, we obtained a certain number of new identities
by fixing $m=1$ in the double-sequences involved. The same
approach was used to study connections between $\sigma_k(x)$
Stirling, $\delta_k(x)$ tanh and $\lambda_(x)$ Lah polynomials.
Finally, we showed that the cumulants of a shifted-gamma
probability density and a negative binomial distribution can be
related by reciprocal linear combinations (Theorem \ref{th2})
which turned out to be an instance of the tanh numbers inversion
formula, hence, that this problem and our problem (Theorem
\ref{conv}) on number arrays, are equivalent.
\bibliographystyle{amsplain}
\providecommand{\bysame}{\leavevmode\hbox
to3em{\hrulefill}\thinspace}
\providecommand{\MR}{\relax\ifhmode\unskip\space\fi MR }
\providecommand{\MRhref}[2]{%
  \href{http://www.ams.org/mathscinet-getitem?mr=#1}{#2}
} \providecommand{\href}[2]{#2}

\end{document}